\documentclass[fleqn]{article}
\usepackage{amssymb,amsmath}
\newtheorem{prop}{Proposition}[section]
\newtheorem{thm}[prop]{Theorem}
\newtheorem{lemma}[prop]{Lemma}
\newtheorem{remark}[prop]{Remark}
\newtheorem{definition}[prop]{Definition}
\newcommand{\vect}{\mbox{vec}}
\newcommand{\emvect}{\mbox{{\em vec}}}
\newcommand{\diag}{\text{diag}}
\newcommand{\emdiag}{\mbox{{\em diag}}}

\begin{document}
\pagestyle{empty}
\Large
\begin{center}
On hidden 
Markov chains and finite stochastic systems
\bigskip \\ \large
Peter Spreij\footnote{Korteweg-de Vries Institute for Mathematics, Universiteit van Amsterdam,
Plantage Muidergracht 24, 1018 TV Amsterdam}
\end{center}
\normalsize
\begin{center}
\today \\
\end{center}
\vspace{5ex}
\begin{abstract}
\noindent
In this paper we study various properties of finite stochastic systems or 
hidden Markov chains as they are alternatively called. We discuss their construction 
following different approaches and
we also derive  recursive filtering 
formulas for the different systems that we consider. The key tool is a simple 
lemma on conditional expectations. \medskip
\end{abstract}

\newpage
~
\newpage

\pagestyle{plain}
\setcounter{page}{1}

\section{Introduction}

In this paper we consider {\em Hidden Markov Chains} (probabilistic functions of a Markov chain) 
that, like the underlying Markov chain, take on finitely many values. The observed process
is denoted by $Y$, the underlying chain by $X$. Hidden 
Markov chains are such that probabilities of future events of $X$ and $Y$ given the 
past only depend on the current state of $X$. Typically this means that $X$ 
satisfies the role of a {\em state process} as it is used in stochastic system 
theory. One of the aims of the present paper is to shed some more light on the 
relation between stochastic systems and hidden Markov chains. There are two slightly 
different
definitions of stochastic systems, related by a time shift of the observed process. 
We will see that a hidden Markov chain satisfies both relations. We will also 
discuss various constructions of a hidden Markov chain. These constructions allow 
different factorization and splitting properties of conditional probabilities of the 
bivariate process $(X,Y)$. We will also study for the different constructions the 
filtering and prediction problems and show that the solutions coincide if one deals 
with a hidden Markov chain in the way we define it. The paper is organized as 
follows. \\
In section~\ref{section:preliminaries} 
we describe the probabilistic behaviour of the joint process $(X,Y)$ in more detail 
using the outer product of $X$ and $Y$ and by using properties of Kronecker products 
of matrices. \\
In section~\ref{section:alternative} we present a somewhat different look at hidden 
Markov chains. It is shown that  certain necessary properties of a hidden Markov 
chain are actually sufficient to construct one. The convenient tool is a simple 
lemma, that is presented in the appendix, 
on conditional expectations that involves a finitely generated $\sigma$-algebra. 
It is also shown that  hidden Markov chains are nothing else but what in 
the engineering literature are called  stochastic systems. 
In particular it is shown that hidden Markov chains satisfy two different notions of 
stochastic systems.
It is also shown how these two 
 notions  are interrelated. This is done in section~\ref{section:stochsyst}. \\
In section~\ref{section:filter} we show how various filtering and prediction 
formulas are simple consequences of the key lemma on conditional expectations of 
the appendix.

\section{Preliminaries} \label{section:preliminaries}

Let $(\Omega,{\cal F},P)$ be a probability space on which all the random 
variables to be encountered below are defined.
Consider the following model for what we will call later a Hidden Markov Chain (HMC).
\begin{eqnarray}
X_t & = & AX_{t-1} + \varepsilon_t, X_0 \label{eqn:system} 
\\
Y_t & = & H_tX_t \label{eqn:observation}  
\end{eqnarray}
Here the {\em state} process $X$ is modelled as a Markov process on the set 
$E=\{e_1,\ldots,e_n\}$ of 
basis vectors of $\mathbb{R}^n$. Moreover, this process is supposed to be 
time-homogeneous with $A$ the matrix of one step transitions probabilities: 
$A_{ij}=P(X_{t+1}=e_i|X_t=e_j)$. The process $\{\varepsilon_t\}$ is then a 
martingale difference sequence adapted to the filtration generated by $X$, 
see \cite[page 17]{elliottetal}. Throughout the paper we assume that each state 
$e_i$ is visited at least once by $X$. If this were not the case, this can always be
accomplished by reducing the state space of $X$ by taking basis vectors of a lower 
dimensional Euclidean space. \\
The {\em observation} or {\em output} process $Y$ takes its values in the set 
$F=\{f_1,\ldots,f_m\}$  
of basis vectors of $\mathbb{R}^m$. The matrices $\{H_t\}$ are assumed to form an {\em iid} sequence, 
independent of $\{X_t\}$, and each column of any of 
these matrices is assumed to be a random element of $F$.  Clearly each $H_t$ 
is the incidence matrix of a random map from $E$ in $F$. Indeed, if 
$Y_t=h_t(X_t)$, with the $h_t$ random maps from $E$ into $F$, then we can 
write $Y_t = \sum_{i=1}^n h_t(e_i)1_{\{X_t = e_i\}}$. So we define 
$H_t=[h_t(e_1),\ldots,h_t(e_n)]$ to get (\ref{eqn:observation}). \\
We will only need the distributions of the colums of $H_t$ (equivalently, the marginal distributions of the $h_t(e_i)$.
These are specified by the expectation $EH_t=G$. We 
assume (without loss of generality) the non-degeneracy condition that none of 
the rows of $G$ is zero.
\\
Define the 
filtration $\mathbb{F} = \{ {\cal F}_t\}$ by ${\cal 
F}_t=\sigma\{X_0,\ldots,X_t,H_0,\ldots,H_t\}$. Clearly both $X$ and $Y$ are 
adapted to this filtration, and so is the sequence $\{\varepsilon_t\}$ which 
is even a martingale difference sequence w.r.t $\mathbb{F}$, because of the 
independence of the sequences $\{X_t\}$ and $\{H_t\}$.\medskip\\
In the current set up, also the joint process $\{(X_t,Y_t)\}$ is Markov. For 
completeness we give its transition probabilities, already given
in~\cite{baumpetrie}, and derive these using simple properties of conditional 
expectations.
\begin{prop}\label{prop:hmm}
The joint process $\{(X_t,Y_t)\}$ is Markov with respect to $\mathbb{F}$ and the 
conditional transition probabilities are given by
\begin{equation}
P(X_{t}=e_i,Y_{t}=f_j|{\cal F}_{t-1})=e_i^\top \mbox{{\em diag}}(AX_{t-1})G^\top f_j 
\label{eq:jointmarkov}
\end{equation}
\end{prop}
{\bf Proof}. Notice first that the indicator of the event 
$\{X_t=e_i,Y_t=f_j\}$ equals $e_i^\top X_tY_t^\top f_j$. Hence we can rewrite the 
conditional probability in equation (\ref{eq:jointmarkov}) as 
$E[e_i^\top X_tY_t^\top f_j|{\cal F}_{t-1}]$. So we compute 
\begin{eqnarray*}
E[X_tY_t^\top |{\cal F}_{t-1}] & = &  E[X_tX_t^\top H_t^\top |{\cal F}_{t-1}] \\
& = & E[E[X_tX_t^\top H_t^\top |{\cal F}_{t-1}\vee\sigma(H_{t})]|{\cal F}_{t-1}] \\
& = & E[E[X_tX_t^\top |{\cal F}_{t-1}\vee\sigma(H_{t})]H_t^\top |{\cal F}_{t-1}] \\
& = & E[E[\mbox{diag}(X_t)|{\cal F}_{t-1}\vee\sigma(H_{t})]H_t^\top |{\cal F}_{t-1}] \\
& = & E[\mbox{diag}(AX_{t-1})H_t^\top |{\cal F}_{t-1}] \\
& = & \mbox{diag}(AX_{t-1})E[H_t^\top |{\cal F}_{t-1}] \\
& = & \mbox{diag}(AX_{t-1})G^\top  
\end{eqnarray*}
The result follows.
\hfill $\square$
\bigskip \\
We will see in section~\ref{section:stochsyst} that it follows from proposition~\ref{prop:hmm} 
that the pair $(X,Y)$ forms a stochastic system in the 
sense of~\cite{picci}. 
\\
We continue with giving an alternative expression for the matrix of one step transition 
probabilities of the joint chain $(X,Y)$. The state 
space of this chain  
consists of all the $nm$ pairs $(e_i,f_j)$. These are renamed and ordered as 
follows: $s_{(j-1)n+i}=(e_i,f_j)$ for $i\in \{1,\ldots,n\}$ and $j\in 
\{1,\ldots,m\}$. Clearly the map $(i,j) \mapsto (i-1)m+j$ is bijective from 
$\{1,\ldots,n\}\times\{1,\ldots,m\}$ onto $\{1,\ldots,nm\}$.       \newline
Instead of working with $(X,Y)$ we will use 
the chain $Z$ that carries the same information and which is  defined by 
$Z_t=\vect(X_tY_t^\top )$. Recall  that
the vec-operator applied to a matrix results in a vector where all 
the columns of this matrix are stacked one underneath the other \cite[p.
30]{magnusneudecker}. Then 
clearly the state space of $Z$ is the set of basis vectors of $\mathbb{R}^{nm}$. 
If we call this set $\{z_1,\ldots,z_{nm}\}$ we see that $(X_t,Y_t)=s_k$ iff 
$Z_t=z_k$. Notice also the following relations. $Z_t=Y_t\otimes X_t$, 
$X_t=({\bf 1}^\top _m\otimes I_n)Z_t$ and $Y_t=(I_m\otimes {\bf 1}^\top _n)Z_t$. Here 
$I_m$ is the $m$-dimensional identity matrix and ${\bf 1}_n$ is the 
$n$-dimensional column vector with all its elements equal to one.
\newline
According to proposition  \ref{prop:hmm} we now get that the $nm \times nm$ 
matrix $Q$  of transition probabilities of $Z$ can be decomposed as a matrix 
with $m^2$ blocks $Q_{ij}$ that are equal to
$\mbox{diag}(G_{i.})A$, where $G_{i.}$ is the $i$-th row of $G$. 
For a more compact formulation we introduce (like in~\cite{spreij}) the following notation. Let
$\Delta(G)$ be the $nm\times n$ matrix defined by
\[ \Delta(G) = \left[
\begin{array}{c}
\diag(G_{1.}) \\
\vdots \\
\diag(G_{m.})
\end{array}\right]             \]
Using the notation $\Delta(G)$ we can now write
\begin{equation}
Q=\Delta(G)A({\bf 1}^\top _m\otimes I_n)               \label{eq:Q}
\end{equation}
In the next lemma we gather some computational results for the 
$\Delta$-operator, that might be of independent interest. Other properties are 
described in~\cite{spreij}.
\begin{lemma}\label{lemma:lemmadelta}
For any matrices $G\in \mathbb{R}^{m\times n}$, 
$M\in \mathbb{R}^{p\times m}$ and $N \in \mathbb{R}^{p\times n}$ and for any 
vectors $w\in \mathbb{R}^n$, $v\in \mathbb{R}^m$  we have
\begin{eqnarray}
MG & = & (M\otimes {\bf 1}^\top _n)\Delta(G)  \label{eq:deltarule} \\
(I_m\otimes \emdiag(w))\emvect(G^\top ) & = & \Delta(G)w \label{eq:rule4}  \\
\emvect(\diag(w)G^\top ) & = & \Delta(G)w \label{eq:rule6}
\end{eqnarray}
\end{lemma}
{\bf Proof.} By direct calculation. \hfill $\square$
\bigskip \\
The expression~(\ref{eq:Q}) for $Q$ can also be obtained through simple matrix manipulations and by application
of lemma~\ref{lemma:lemmadelta}. By definition of $Q$ we have 
$E[Z_{t+1}|{\cal F}_t]=QZ_t$. So we compute the conditional expectation
\begin{eqnarray*}
E[Z_{t+1}|{\cal F}_t] & = & E[\vect(X_{t+1}Y_{t+1}^\top )|{\cal F}_t] \\
& = & \vect(E[X_{t+1}Y_{t+1}^\top |{\cal F}_t]) \\
& = & \vect(E[X_{t+1}X_{t+1}^\top H^\top _{t+1}|{\cal F}_t]) \\
& = & \vect(\diag(AX_t)G^\top ) \\
& = & (I_m\otimes \diag(AX_t))\vect(G^\top ) \\
& = & \Delta(G)AX_t \\
& = & \Delta(G)A({\bf 1}^\top _m\otimes I_n)Z_t  \\
& = & QZ_t
\end{eqnarray*}
Here we used in the fifth equality a known result for the vec-operator of 
the product of three matrices (see 
\cite[page 30]{magnusneudecker}) and in the sixth equality equation~(\ref{eq:rule4}). 
\medskip \\
If the vector $p_0$ represents the initial distribution of $X$, then 
the initial distribution of $Z$ is given by the vector 
$EZ_0=\mbox{vec}(\diag(p_0)G^\top )$: $EZ_0=E\vect(X_0Y_0^\top ) = 
\vect(E\diag(X_0)H_0^\top ) = \vect(\diag(p_0)G^\top )$, since $X_0$ and $H_0$ are independent. Notice that 
$\vect(\diag(p_0)G^\top )=\Delta(G)p_0$, because of (\ref{eq:rule6}). \newline
Similarly one can show that $\Delta(G)\pi$ is an invariant probability vector for $Z$, 
if $X$ has an invariant probability vector $\pi$.
\medskip \\ 
It is easy to see from equation~(\ref{eq:jointmarkov}) that the 
"{\em factorization property}"~\cite{finesso} holds:
\begin{equation}
P(X_{t}=e_i,Y_{t}=f_j|{\cal F}_{t-1})=P(Y_t=f_j|X_t=e_i) P(X_t=e_i|X_{t-1}) 
\label{eq:factor}
\end{equation}
The proof is straight forward from proposition~\ref{prop:hmm} (used in the first equality below):
\begin{eqnarray*} 
P(X_{t}=e_i,Y_{t}=f_j|{\cal F}_{t-1}) & = & 
e_i^\top \mbox{diag}(AX_{t-1})G^\top f_j \\
& = & (AX_{t-
1})^\top \mbox{diag}(e_i)G^\top f_j \\ 
& = & (AX_{t-
1})^\top e_ie_i^\top G^\top f_j  \\
& = & 
P(X_t=e_i|X_{t-1})G_{ji} \\
& = & P(X_t=e_i|X_{t-
1})P(Y_t=f_j|X_t=e_i).
\end{eqnarray*}
Using the matrix $\Delta(G)$ and lemma \ref{lemma:lemmadelta} we can also compactly rephrase the factorization 
property (\ref{eq:factor}). It becomes
\begin{equation}
E[Z_t|{\cal F}_{t-1}] = \Delta(G)E[X_t|X_{t-1}], \forall t.
\label{eq:factor2}
\end{equation}
This can be verified as follows. First, 
using proposition~\ref{prop:hmm} again, we rewrite (\ref{eq:factor}) as 
\[ 
e_i^\top E[X_tY_t^\top|\mathcal{F}_{t-1}]f_j = G_{ji}e_i^\top 
E[X_t|\mathcal{F}_{t-1}].\]
Since the right hand side of this equality equals
$f_j^\top G \diag(e_i)E[X_t|\mathcal{F}_{t-1}]$, which is equal to
$f_j^\top G \diag(E[X_t|\mathcal{F}_{t-1}])e_i$,  we get 
\[ 
E[X_tY_t^\top |{\cal F}_{t-1}] = \diag(E[X_t|X_{t-1}])G^\top .
\] 
Since 
$\vect(X_tY_t^\top ) = Z_t$ and 
\[
\vect(\diag(E[X_t|X_{t-
1}])G^\top )=(I_m\otimes\diag(E[X_t|X_{t-1}]))\vect(G^\top ),
\]
 we use (\ref{eq:rule4}) 
to write the RHS of this last equality as $\Delta(G)E[X_t|X_{t-1}]$, from which~(\ref{eq:factor2}) follows. 
\begin{remark}
\textup{
 The validity of equation (\ref{eq:factor2}) has been seen to be a consequence 
of the special form of the transition matrix $Q$ in (\ref{eq:Q}). But also the 
converse holds. If (\ref{eq:factor2}) holds, we get at once that $Z$ is 
$\mathbb{F}$-Markov, if $X$ is $\mathbb{F}$-Markov. And if we denote the 
transition matrix of $Z$ by $Q$ and that of $X$ by $A$, we automatically get 
$(\ref{eq:Q})$ back. See proposition~\ref{prop:properties}.
}
\end{remark} 
As an alternative to looking at the bivariate process $(X,Y)$ via the process $Z$ as above, we 
study the process $W$, again built from $X$ and $Y$ and defined by $W_t=Y_{t-1}\otimes 
X_t$ for $t \geq 1$. Along with this process we consider the filtration $\mathbb{G}$ of 
$\sigma$-algebras  
${\cal G}_t := \sigma\{H_0,\ldots,H_{t-1},X_0,\ldots,X_t\}$. Then $W$ is 
$\mathbb{G}$-adapted and the ${\cal G}_t$ and the ${\cal F}_t$ are related by
${\cal F}_{t-1}\vee\sigma(X_t) = {\cal G}_t$ and ${\cal G}_t\vee\sigma(H_t)= 
{\cal F}_t$. \\
Then
by similar computations as we carried out before and by using the Markov 
property  of $Z$ we obtain the relations
\begin{eqnarray}
E[W_t|{\cal F}_{t-1}] & = & (I_m\otimes A)Z_{t-1}, \label{eq:zandw} \\
E[W_t|{\cal G}_{t-1}] & = & GX_{t-1} \otimes AX_{t-1} = (I_m\otimes A)\Delta(G) 
X_{t-1}.  
\label{eq:gmarkov}
\end{eqnarray}
In particular it follows that $W$ is $\mathbb{G}$-Markov (and hence the pair $(X,Y)$
is a stochastic 
system in the sense of~\cite{schuppen}, see 
section~\ref{section:stochsyst}) with transition matrix 
\begin{equation}
R:= (I_m\otimes A)\Delta(G) ({\bf 1}^\top _m\otimes I_n).
\label{eq:r}
\end{equation}
Observe also that $W$ has the {\em splitting property}
\begin{equation}
E[W_{t+1}|{\cal G}_t] = E[Y_t|{\cal G}_t]\otimes E[X_{t+1}|{\cal G}_t],
\label{eq:splitting}
\end{equation}
which immediately follows from (\ref{eq:gmarkov}). 
\begin{remark}
\textup{
 The assumption in this section that the sequence $\{H_t\}$ is 
{\em iid} with $EH_t=G$ can in  
principle be relaxed to assuming that $\{H_t-G\}$ is a martingale difference
sequence with respect to its own filtration without changing the results of this section. However, this only 
{\em appears} to be a relaxation, in fact they are equivalent assumptions. Indeed, 
let $\{H_t-G\}$ be a martingale difference sequence and consider 
$k_t=\vect(H_t)$. Then $k_t$ takes its values in the set of basis vectors of 
$\mathbb{R}^{mn}$ and $\{k_t-\vect(G)\}$ is again a martingale difference 
sequence. Let $e$ be one of these basis vectors. Then  
$P(k_{t+1}=e|k_0,\ldots,k_t)=e^\top E[k_{t+1}|k_0,\ldots,k_t]=e^\top \vect(G)$, 
which doesn't depend on $k_0,\ldots,k_t$, nor on time. Hence $\{k_t\}$ is an iid sequence 
and so is $\{H_t\}$. \\
We can also replace
(\ref{eqn:observation}) with the equivalent equation
\begin{equation}\label{eq:y'}
Y_t=GX_t + \eta_t
\end{equation}
where $\eta$ forms a martingale difference sequence with respect to $\{{\cal 
F}_t\}$, and it even holds that $\eta_t=Y_t-E[Y_t|\sigma(X_t)\vee{\cal F}_{t-1}]=Y_t-E[Y_t|{\cal G}_{t}]$. The combined
set of equations~(\ref{eqn:system}) and (\ref{eq:y'}) are of the form that is 
commonly used in (stochastic) systems theory.
We will come back to stochastic systems in section~\ref{section:stochsyst}.
}
\end{remark}
\begin{remark}\label{remark:law} 
\textup{
As a final remark we notice that all the properties mentioned 
above in terms of conditional expectations given the $\sigma$-algebras ${\cal 
F}_t$ and ${\cal G}_t$ remain valid if we replace the former one with 
$\sigma\{X_0,\ldots,X_t,Y_0,\ldots,Y_t\}$ and the latter one with 
$\sigma\{X_0,\ldots,X_t,Y_0,\ldots,Y_{t-1}\}$. Hence the law of the bivariate 
process $(X,Y)$, being a Markov chain with respect to its own filtration, is 
completely specified by the matrices $A$ and $G$ and the initial law of $X$. 
It follows that any bivariate Markov process $(X,Y)$, that is such that the 
transition matrix $Q$ of the associated process $Z=Y\otimes X$ is of the form 
$(\ref{eq:Q})$ and that has initial law $EZ_0=\Delta(G)p_0$ where $p_0=EX_0$, 
can be constructed as the output of the system  
(\ref{eqn:system}) and (\ref{eqn:observation}).
}
\end{remark}
In view of remark~\ref{remark:law} above we adopt the following
\begin{definition}\label{def:hmc}
A  bivariate process $(X,Y)$ that assumes finitely many values is called a Hidden Markov Chain (HMC) if the process 
$Z=Y\otimes X$ is Markov with respect to the filtration 
$\mathbb{F}=\{\mathcal{F}_t\}$ defined by
$\mathcal{F}_t=\sigma\{X_0,\ldots,X_t,Y_0,\ldots,Y_t\}$ and if its matrix of 
transition probabilities is given by~(\ref{eq:Q}).
\end{definition}

\section{Alternative descriptions of a HMC} 
\label{section:alternative} 

There are various ways to describe some properties of a stochastic system or a 
Hidden Markov chain. We mention a few possibilities and show how these can be used as 
building stones for a HMC. \\ Let $X$ and $Y$ be two stochastic processes taking 
values in the sets $E$ and $F$ respectively, like in section~\ref{section:preliminaries}. Let $Z$  again be 
the process $Y\otimes X$. For the time being no further assumptions on $X$ and $Y$ 
are imposed, except that redundant states are excluded in the sense that each state of $X$ is visited at least 
once with probability one and likewise for $Y$. \\
In this section (and all subsequent ones) we  assume that for all $t$ the $\sigma$-algebra 
${\cal F}_t$ is generated by $X_0,\ldots,X_t,Y_0,\ldots,Y_t$. The family $\{ {\cal 
F}_t\}$ is again denoted by $\mathbb{F}$. We also consider the process $W$ again, 
with $W_t=Y_{t- 1}\otimes X_t$, adapted to the filtration $\mathbb{G}=\{{\cal 
G}_t\}$, with ${\cal G}_t$ generated by $X_0,\ldots,X_t,Y_0,\ldots,Y_{t-1}$. Notice 
again the relations
\begin{eqnarray*}
\mathcal{F}_t & = & \mathcal{G}_t\vee\sigma(Y_t) \\
\mathcal{G}_t & = & \mathcal{F}_{t-1}\vee\sigma(X_t) \\
\end{eqnarray*}

\subsection{Alternative description of Z} \label{subsection:altz}

We now list a set possible properties that the processes $X$, $Y$ and $Z$ may 
possess. 
\begin{enumerate} 
\item\label{item:zxmarkov}
The process $Z$ is time homogeneous $\mathbb{F}$-Markov with matrix $Q$ of transition 
probabilities, so $E[Z_{t+1}|{\cal F}_t]=QZ_t$. Moreover we assume that this
conditional expectation only depends on $X_t$, which implies that there
exists a matrix $ \bar{Q}$ such that $E[Z_{t+1}|{\cal F}_t]=\bar{Q}X_t=\bar{Q}({\bf 1}^\top _m\otimes 
I_n)Z_t, \forall t$. Hence $Q=\bar{Q}({\bf 1}^\top _m\otimes 
I_n)$. 
\item\label{item:outp} 
The {\em output property} holds:  
\[
E[Y_t|\mathcal{G}_t] = E[Y_t|{\cal F}_{t-1}\vee 
\sigma(X_t)] =E[Y_t|\sigma(X_t)], \forall t.
\]
If this property holds, we use the matrix $G$ 
defined by                                     
$E[Y_t|\sigma(X_t)]=GX_t$, where we also assume that $G$ is not depending on $t$.  
$G$ is then such that the columns $G_{\cdot i}$ are equal to
$E[Y_t|X_t=e_i]$.
\item\label{item:eoutp}
The {\em extended output property} holds: 
\begin{equation}\label{eq:eoutp}
E[Z_t|\mathcal{G}_t] = E[Z_t|{\cal F}_{t-1}\vee  
\sigma(X_t)] = E[Z_t|\sigma(X_t)], \forall t.
\end{equation} 
In this case we define the matrix $B$ (assumed to be independent of $t$)
by                                     
$E[Z_t|\sigma(X_t)]=BX_t$. 
\item\label{item:factor}
The {\em factorization property} holds: There exists a matrix $K \in 
\mathbb{R}^{m\times n}$ such that
\begin{equation}
E[Z_t|{\cal F}_{t-1}] = \Delta(K) E[X_t|{\cal F}_{t-1}], \forall t 
\label{eq:condz}
\end{equation}
\end{enumerate}
First we comment on the factorization property. We showed that it is valid for the 
HMC of section~\ref{section:preliminaries}.But 
one can always factorize  $E[Z_t|{\cal F}_{t-1}]$ with a second
factor $E[X_t|{\cal F}_{t-1}]$ as in (\ref{eq:condz}), however in general the left factor is a random 
(${\cal F}_{t-1}$-measurable) diagonal matrix, see equation~(\ref{eq:roughfactor}) below.
\\
Denote by $P_i$ the conditional measure on $(\Omega, {\cal F})$ given 
$X_t=e_i$.  Expectation with respect 
to these measures will be denoted by $E_i$, with the understanding that expectations $E_iU$ are set equal to zero, if
$P(X_t=e_i)=0$ (cf.\ the appendix).
Then for any sub-$\sigma$-algebra ${\cal F}^0$ of ${\cal F}$ and any 
integrable random variable $U$ we have from equation (\ref{eq:bayes2}) in the appendix the relation
\begin{equation}
E[U1_{\{X_t=e_i\}}|{\cal F}^0] = E_i[U|{\cal F}^0] P(X_t=e_i|{\cal F}^0).
\label{eq:condT}
\end{equation}
Application  of equation (\ref{eq:condT}) with $U=Y_t^\top $, ${\cal F}^0 = {\cal 
F}_{t-1}$ for all $i$ yields
\begin{equation}
E[X_tY_t^\top |{\cal F}_{t-1}] = \diag(E[X_t|{\cal F}_{t-1}]) E^\top ,
\end{equation}
where  $E^\top $ is the transpose of the matrix $E$ that has columns $E_i[Y_t|{\cal F}_{t-
1}]$. Apply then (\ref{eq:rule6}) to get
\begin{equation}
E[Z_t|{\cal F}_{t-1}]=\Delta(E) E[X_t|{\cal F}_{t-1}]. 
\label{eq:roughfactor}
\end{equation}

\begin{prop} 
Properties \ref{item:outp}, \ref{item:eoutp} and \ref{item:factor} are equivalent. 
Moreover the matrices $B$, $G$ and $K$ are 
related via $B=\Delta(G)$ and $K=G$. 
\label{prop:properties}
\end{prop}
{\bf Proof.} 
Trivially the output property \ref{item:outp} follows from the extended output property~\ref{item:eoutp} 
by left multiplication with $I_m\otimes {\bf 1}^\top _n$.  \\
Conversely, assume that the output property holds. Then we have $E[Z_t|{\cal F}_{t-
1}\vee\sigma(X_t)] = E[Y_t|{\cal F}_{t-1} \vee \sigma(X_t)]\otimes X_t = 
E[Y_t|X_t]\otimes X_t= E[Z_t|X_t]$, which shows that the extended output property holds. \\
To see the relation between $B$ and $G$, notice that in this case we have 
$BX_t=E[Z_t|X_t]=E[Y_t|X_t]\otimes X_t = GX_t\otimes X_t = \vect(X_tX_t^\top G^\top ) 
= \vect(\diag(X_t)G^\top ) = \Delta(G)X_t$. Here we used the usual relations 
between the vec-operator and Kronecker products as well as (\ref{eq:rule6}) 
in the last equality. \\
Assume that the extended output property holds. Use then reconditioning 
in~(\ref{eq:eoutp}) to get:
$E[Z_t|{\cal F}_{t-1}] = 
E[E[Z_t|{\cal F}_{t-1}\vee \sigma(X_t)]|{\cal F}_{t-1}] = E[E[Z_t|X_t]|{\cal 
F}_{t-1}] = BE[X_t|{\cal F}_{t-1}]$. It follows from~(\ref{eq:roughfactor}) that $B=\Delta(E)$, but
since $B$ is nonrandom, the validity of the factorization property follows. \\
Conversely, assume that the factorization property~\ref{item:factor} holds. 
Take expectations in (\ref{eq:condz}). Then $EZ_t=\Delta(K)EX_t$. From 
the definition of $G$ (in property~\ref{item:outp}) we get 
$EZ_t=EE[Y_t|\sigma(X_t)]\otimes X_t]=E (GX_t\otimes X_t)=\Delta(G)EX_t$. 
Since for each $i$ there is a $t$ such that the $i$-th component of $EX_t$ is strictly positive,  
it follows from the blockwise diagonal structure of the 
$\Delta$-matrices that $\Delta(G)=\Delta(K)$ and $G=K$. 
\\
Next we show that the output property holds. Assume for a moment that all elements of 
$EX_t$ are positive. According to equations~(\ref{eq:bayes}) and~(\ref{eq:condi}) we 
have
\[
E[Y_t|\mathcal{F}_{t-1}\vee\sigma(X_t)] = \sum_i\frac{E[Y_te_i^\top X_t|\mathcal{F}_{t-1}]}{e_i^\top 
E[X_t|\mathcal{F}_{t-1}]}e_i^\top X_t.
\]
Since $Y_te_i^\top X_t= (I_m\otimes e_i^\top)Z_t$ and using the factorization 
property, we can rewrite this as
\[
\sum_i \frac{(I_m\otimes e_i^\top)\Delta(G)E[X_t|\mathcal{F}_{t-1}]}{e_i^\top 
E[X_t|\mathcal{F}_{t-1}]}e_i^\top X_t.
\]
Because $(I_m\otimes e_i^\top)\Delta(G)=Ge_ie_i^\top $, this reduces to
\[
G\sum_i \frac{e_ie_i^\top E[X_t|\mathcal{F}_{t-1}]}{e_i^\top 
E[X_t|\mathcal{F}_{t-1}]}e_i^\top X_t,
\]
which in turn is nothing else but $GX_t$, from which we obtain the output property.
In the case where the vector $EX_t$  has some elements equal to zero, the above 
procedure is still valid, provided we let the summation indices  run through the set 
$\{i:e_i^\top EX_t >0\}$.
\hfill$\square$
\medskip \\
Similar to what we found in the previous section we have
\begin{prop}\label{prop:zxmarkov}
Assume that the factorization property
\ref{item:factor} holds (or, equivalently in view of proposition~\ref{prop:properties}, 
the output or extended output property). 
Then the  following two statements are 
equivalent. \\ 
(i) $Z$ is $\mathbb{F}$-Markov with transition matrix $Q=\bar{Q}({\bf 1}^\top _m\otimes 
I_n)$. \\
(ii) $X$ is $\mathbb{F}$-Markov with transition matrix $A$  \\ 
Furthermore we have in each of these situations the relation 
$\bar{Q}=\Delta(G)A$. 
\end{prop}
{\bf Proof.}  (i) $\Rightarrow$ (ii): Clearly $X$ is  
$\mathbb{F}$-Markov with transition matrix $A=({\bf 1}^\top _m\otimes I_n)\bar{Q}$ and 
then it follows from the factorization property
that $\bar{Q}X_{t-1}= E[Z_t|{\cal F}_{t-1}] = \Delta(G) E[X_t|{\cal F}_{t-
1}]= \Delta(G)AX_{t-1}$. \\
Conversely, (ii) $\Rightarrow$ (i) follows in a similar way. $E[Z_t|{\cal 
F}_{t-1}] = \Delta(G)E[X_t|{\cal F}_{t-1}]= \Delta(G)AX_{t-1} = 
\Delta(G)A({\bf 1}^\top _m\otimes I_n)Z_{t-1}$, so $Z$ is $\mathbb{F}$-Markov with 
transition matrix $Q=\Delta(G)A({\bf 1}_m^\top \otimes  I_n)$.
\hfill $\square$
\begin{remark}
\textup{
The main implication of proposition~\ref{prop:zxmarkov} is that the proces $Z$
is a Markov chain whose transition probabilities only depend on the past value of $X$, 
if one starts out with a $\mathbb{F}$-Markov chain $X$ and imposes that the output condition 
holds. Clearly, if $X$ is just Markov with respect to its own filtration and if the 
factorization property is replaced with the stronger condition 
$E[Z_t|{\cal F}_{t-1}] = \Delta(K) E[X_t|{\cal F}_{t-1}]$, the same conclusion 
follows.
}
\end{remark}
\begin{remark}\label{remark:zwmarkov}
\textup{
We also observe, like in section \ref{section:preliminaries}, that 
the fact that $Z$ is $\mathbb{F}$-Markov with $\bar{Q}=\Delta(G)A$ implies that 
$W$ is $\mathbb{G}$-Markov, with transition matrix $R=(I_m\otimes A)\Delta(G)$. 
One easily checks that with the present choice of the filtrations equations 
(\ref{eq:zandw}) and (\ref{eq:gmarkov}) remain valid, and that in particular 
the  factorization property holds.
}
\end{remark}

\subsection{Alternative description of W} \label{subsection:altw}

Like in subsection \ref{subsection:altw}, we can also list a set of desirable 
properties of $W$. Consider thereto
\begin{enumerate}
\item
$W$ is a time homogeneous $\mathbb{G}$-Markov chain with a transition matrix $R$. Moreover, we have 
that conditional expectation $E[W_{t+1}|\mathcal{G}_t]$ depends only on $X_t$. This 
means that there is a matrix $\bar{R}$ such that $R=\bar{R}({\bf 1}^\top_m\otimes 
I_n)$.
\item
The {\em splitting property} holds: 
\begin{equation}
E[W_t|{\cal 
G}_{t-1}] = E[Y_{t-1}|{\cal G}_{t-1}]\otimes E[X_t|{\cal G}_{t-1}], \forall 
t. \label{eq:splitprop}
\end{equation}
\end{enumerate}
Then we have similar to proposition \ref{prop:zxmarkov}
\begin{prop}
Under the splitting property (\ref{eq:splitprop}) there is equivalence 
between \\
(i) $W$ is $\mathbb{G}$-Markov with transition matrix $R=\bar{R}({\bf
1}_m^\top \otimes I_n)$. \\
(ii) $X$ is $\mathbb{G}$-Markov with a transition matrix $A$. \\
Moreover, in each of these cases we have the relation $\bar{R}=(I_m\otimes 
A)\Delta(G)$.
\label{prop:wxmarkov}
\end{prop}
{\bf Proof.} We omit the proof of proposition \ref{prop:wxmarkov}, since it is similar to 
that of proposition \ref{prop:zxmarkov}.
\begin{remark}
\textup{
The main message of proposition~\ref{prop:wxmarkov} is that to have $W$ Markov with transition probabilities only 
depending on past values of $X$ it is sufficient to start with a $\mathbb{G}$-Markov chain $X$ 
and to assume that the splitting property~(\ref{eq:splitprop}) holds.
}
\end{remark}
\begin{remark}\label{remark:wzmarkov}
\textup{
We noticed in remark~\ref{remark:zwmarkov}, that from the assumption that $Z$ 
is $\mathbb{F}$-Markov and the validity of the factorization property, one could deduce that 
$W$ is $\mathbb{G}$-Markov. Conversely, given that 
$W$ is $\mathbb{G}$-Markov with the transition matrix as in (\ref{eq:r}) above, 
we can also  
deduce that $Z$ is $\mathbb{F}$-Markov with $Q$ as in (\ref{eq:Q}) as its 
transition matrix (and that  
equation (\ref{eq:factor2}) holds). This also follows from more general 
considerations to be explained at the end of  section 
\ref{section:stochsyst}, but here we give an explicit calculation.
\\
So let  $W$ be a 
$\mathbb{G}$-Markov process with transition matrix $R=\bar{R}({\bf 1}_m^\top  \otimes 
I_n)$. Then 
$E[W_t|\mathcal{G}_{t-1}]=  \bar{R}X_{t-1}$.
From this it follows that
\[
E[X_t|\mathcal{G}_{t-1}]=({\bf 1}_m^\top  \otimes I_n)\bar{R}X_{t-1}= AX_{t-1}
\]
with $A= ({\bf 1}_m^\top  \otimes I_n)\bar{R}$. Furthermore we have
$E[Y_t|\mathcal{G}_t]=(I_m\otimes {\bf 1}_n^\top )E[ W_{t+1}|\mathcal{G}_t] = 
(I_m\otimes {\bf 1}_n^\top )\bar{R} X_t= GX_t$ with $G=(I_m\otimes {\bf 1}_n^\top )\bar{R}$. \\
We now  compute $E[Z_{t+1}|{\cal F}_t] =E[E[Y_{t+1}|{\cal 
G}_{t+1}]\otimes X_{t+1}|{\cal F}_t]$. By the relation that we just showed, this becomes 
$E[GX_{t+1}\otimes X_{t+1}|{\cal F}_t]$ 
which is $\Delta(G)E[X_{t+1}|{\cal F}_t]$. We have reached our goal as soon as we show 
that $E[X_{t+1}|{\cal F}_t]=E[X_{t+1}|{\cal G}_t]$.
But it is easy to see that this follows immediately from
the splitting property (actually it is equivalent). 
\\ 
Thus we showed the Markov 
property of $Z$ with respect to $\mathbb{F}$ and found its 
transition matrix.
}
\end{remark}
Altogether we summarize our findings of this section in
\begin{thm}\label{thm:summary}
There is equivalence between\\
(a) $X$ is $\mathbb{F}$-Markov and the factorization property holds.  \\
(b) $X$ is $\mathbb{G}$-Markov and the splitting property holds. \\
(c) $Z$ is $\mathbb{F}$-Markov with transition matrix $Q$ as in (\ref{eq:Q}). \\
(d) $W$ is $\mathbb{G}$-Markov with transition matrix $R$ as in (\ref{eq:r}). \\
(e) $(X,Y)$ is a hidden Markov chain.
\end{thm}
{\bf Proof.} The equivalence of (a) and (c) is just 
proposition~\ref{prop:zxmarkov}, that of (b) and (d) is proposition~\ref{prop:wxmarkov}. 
Equivalence of (c) and (d) is the content of remarks~\ref{remark:zwmarkov} 
and~\ref{remark:wzmarkov}, whereas (c) and (e) are equivalent by definition~\ref{def:hmc} of a hidden Markov chain.
\hfill$\square$

\section{Stochastic systems} \label{section:stochsyst}

In the previous sections we restricted ourselves  to  time 
homogeneous processes, implying that all conditional probabilities and expectations don't depend on time directly. 
In the present section where explicit calculations are absent, this restriction playes no role. 
We introduce some notation.
Given a stochastic process $\zeta$ with values in some arbitrary measurable 
space, we denote for all $t$ by ${\cal F}^\zeta_t$ the 
$\sigma$-algebra generated by the $\zeta_s$ for $s\leq t$ and by ${\cal 
F}^{\zeta+}_t$ the $\sigma$-algebra generated by the $\zeta_s$ for $s \geq t$. Many 
of the results in the previous sections can be abstractly formulated in terms of 
properties of {\em stochastic systems}.
A stochastic system is a formally defined concept. The main ingredients are a {\em state 
process} $X$ and an {\em output process} $Y$ (defined on a suitable 
probability space and taking values in some other spaces) and certain conditional independence relations. \\
Let us therefore recall some facts on conditional independence. 
Two $\sigma$-algebras ${\cal H}_1$ and ${\cal H}_2$ are called conditionally 
independent given a $\sigma$-algebra ${\cal G}$ if for all bounded ${\cal 
H}_i$-measurable functions $H_i$ ($i=1,2$) the relation $E[H_1H_2|{\cal G}] = 
E[H_1|{\cal G}]E[H_2|{\cal G}]$ holds. A convenient characterization of this 
is that $\sigma$-algebras ${\cal H}_1$ and ${\cal H}_2$ are  conditionally 
independent given $\sigma$-algebra ${\cal G}$ if for all bounded ${\cal 
H}_1$-measurable functions $H_1$  the relation $E[H_1|{\cal G}\vee {\cal H}_2] 
= 
E[H_1|{\cal G}]$ holds.\\
In the literature one can 
find two definitions of a stochastic system, that are slightly different. The first one is due to 
Picci \cite{picci}, 
and the essential part of the definition is that for all $t$ the 
$\sigma$-algebras 
${\cal F}^{X+}_{t+1}\vee {\cal F}^{Y+}_{t+1}$ and ${\cal F}^{X}_{t}\vee {\cal 
F}^{Y}_{t}$ are conditionally independent given $\sigma(X_t)$. The other one 
is due to Van Schuppen \cite{schuppen} 
in which the conditional independence relation between 
$\sigma$-algebras becomes: for all $t$ the $\sigma$-algebras ${\cal 
F}^{X+}_{t}\vee {\cal F}^{Y+}_{t}$ and  
${\cal F}^{X}_{t-1}\vee {\cal 
F}^{Y}_{t-1}$ are conditionally independent given $\sigma(X_t)$. 
Implications of the two different definitions for the filtering problem will 
be discussed in section~\ref{section:filter}. \\
We will write $(X,Y) \in \Sigma_P$ if the pair of processes $(X,Y)$ is a 
stochastic system according to~\cite{picci} and $(X,Y) \in 
\Sigma_S$ if it is one in the sense of~\cite{schuppen}. 
Using this notation, we see that $(X,Y)\in\Sigma_P$ is equivalent with saying 
that $Z$ is an $\mathbb{F}$-Markov process with transition probabilities 
depending on $X$ only, and 
that $(X,Y)\in\Sigma_S$ is equivalent with saying that 
$W$ is a $\mathbb{G}$-Markov process with transition 
probabilities depending on $X$ only.
Notice that both for a stochastic system $(X,Y)$ either in $\Sigma_P$ or in $\Sigma_S$
the state process is always Markov relative to its own 
filtration.    
\medskip \\
An obvious relation between the different concepts is that 
$(X,Y) \in \Sigma_P$ iff  $(X,\sigma Y) \in \Sigma_S$, where $\sigma Y$ is 
the process defined by  
$\sigma Y_t=Y_{t+1}$. Another relation is given in the following
\begin{prop}
A pair $(X,Y)$ belongs to $\Sigma_S$ and the splitting property holds 
iff it 
belongs to $\Sigma_P$ and the output property (or the factorization 
property) 
holds. 
\end{prop}
{\bf Proof.} Suppose that $(X,Y) \in \Sigma_P$ and that the output property holds. 
Since $X$ is $\mathbb{F}$-Markov, we have $E[X_{t+1}|{\cal F}_t] = 
E[X_{t+1}|X_t]$, which is ${\cal G}_t$ measurable and therefore equal to 
$E[X_{t+1}|{\cal G}_t]$, which is equivalent to the splitting property because of the characterization
of conditional independence given at the beginning of this section. \\
Next we show that $(X,Y)$ also belongs to $\Sigma_S$. We compute
\begin{eqnarray*}
E[W_{t+1}|{\cal G}_t] & = & E[E[W_{t+1}|{\cal F}_t]|{\cal G}_t] \\
& = & E[Y_t\otimes E[X_{t+1}|{\cal F}_t]|{\cal G}_t] \\
& = &
  E[Y_t\otimes E[X_{t+1}|X_t]|{\cal G}_t] \\
  & = & E[Y_t|{\cal G}_t] \otimes
E[X_{t+1}|X_t],
\end{eqnarray*}
which is $\sigma(X_t)$-measurable, because of the output property. \\
Conversely, letting $(X,Y) \in \Sigma_S$ we automatically get the output 
property, because $E[Y_t|\mathcal{G}_t]=(I_n\otimes{\bf 1}^\top_n)E[W_{t+1}|\mathcal{G}_t] 
=(I_n\otimes{\bf 1}^\top_n)E[W_{t+1}|X_t]$ in view of $(X,Y)\in\Sigma_S$. 
Assuming the conditional independence relation we obtain the 
Markov property of $Z$ from 
\begin{eqnarray*}
E[Z_{t+1}|{\cal F}_t] &  = & E[E[Z_{t+1}|{\cal 
G}_{t+1}]|{\cal F}_t] \\
&  = & E[E[Y_{t+1}|{\cal G}_{t+1}]\otimes X_{t+1}|{\cal 
F}_t] \\
& = &  E[E[Y_{t+1}|X_{t+1}]\otimes X_{t+1}|{\cal F}_t] \text{ (output property)}\\
& = &
E[E[Y_{t+1}|X_{t+1}]\otimes X_{t+1}|{\cal G}_t] \text{ (splitting property)}\\
&  = &
E[E[Y_{t+1}|X_{t+1}]\otimes X_{t+1}|X_t] \text{ ($W$ is $\mathbb{G}$-Markov)},
\end{eqnarray*} which shows that $(X,Y) \in 
\Sigma_P$. \hfill $\square$
\begin{remark}
\textup{
 Observe that we already encountered a computational form of this 
proposition in subsections \ref{subsection:altz} and  \ref{subsection:altw}.
}
\end{remark}
The connection between systems in $\Sigma_P$ and $\Sigma_S$ and Hidden Markov chains 
is described as
\begin{prop}
A finite valued time homogeneous system belonging both to $\Sigma_P$ and to $\Sigma_S$ is a Hidden Markov chain and 
vice versa.
\end{prop}
{\bf Proof.}
If $(X,Y)$ is a HMC, then it follows from theorem~\ref{thm:summary} that it belongs 
to both $\Sigma_P$ and $\Sigma_S$. The converse statement follows in a similar way from this 
theorem.
\hfill$\square$

\section{Filtering} \label{section:filter}
In this section we give some filtering and prediction formulas. By the 
filtering problem for a system $(X,Y)$ belonging to $\Sigma_P$ or to $\Sigma_S$ 
we mean the determination for each $t$ of the conditional law of $X_t$ 
given $Y_0,\ldots,Y_t$. As before, for each $t$ we denote by ${\cal F}^Y_t$ the 
$\sigma$-algebra generated by $Y_0,\ldots,Y_t$. Since the state space of $X$ is a set of basis 
vectors, this conditional law is completely determined by the conditional 
expectation $E[X_t|{\cal F}^Y_t]$. The prediction problem is to determine for 
each $t$ the conditional law of $X_{t+1}$ given $Y_0,\ldots,Y_t$, that is 
completely characterized by the conditional expectations $E[X_{t+1}|{\cal 
F}^Y_t]$. We will use the notations $E[X_t|{\cal F}^Y_t]=\hat{X}_t$ and 
$E[X_{t+1}|{\cal  F}^Y_t]= \hat{X}_{t+1|t}$. Similarly we write 
$E[Y_{t+1}|{\cal  F}^Y_t]= \hat{Y}_{t+1|t}$. In addition to the above one wants to 
have $\hat{X}_t$ and $\hat{X}_{t+1|t}$ in recursive form. We shall see below that 
the recursions for the cases $(X,Y)\in \Sigma_P$ and $(X,Y)\in \Sigma_S$  are different.
\\
In the book \cite{elliottetal} recursive formulae for  
unnormalized filters are obtained by a measure transformation. Here we 
undertake a direct approach, that leads to a simple recursive formula for the 
conditional probabilities itself. The key argument is in all cases provided by 
lemma~\ref{lemma:bayes}.

\subsection{Filter for $\Sigma_P$} \label{subsection:filterp}
In this section we obtain the filter for a system in $\Sigma_P$, so we work 
with a Markov chain $Z_t=X_t\otimes Y_t$ with transition matrix 
$Q=\bar{Q}({\bf 1}_m^\top \otimes I_n)$. The matrix $\bar{Q}$
we can write as
\begin{equation}
\bar{Q}=
\left[
\begin{array}{c}
Q_1 \\
\vdots \\
Q_m
\end{array}
\right].
\end{equation}
with the $Q_i$  in $\mathbb{R}^{n\times n}$. No further assumptions 
on the $Q_i$ are made. Observe that the $Q_i$ have the interpretation that 
\begin{equation}\label{eq:intQ}
Q_iX_t =  E[X_{t+1}1_{\{Y_{t+1}=f_i\}}|{\cal F}_t].
\end{equation}
We have the following result (alternatively presented in \cite{picci}).
\begin{thm}
The filter $\hat{X}$ is given by the recursion
\begin{equation}
\hat{X}_t=
\left[
\begin{array}{ccc}
\frac{Q_1\hat{X}_{t-1}}{{\bf 1}_n^\top Q_1\hat{X}_{t-1}} & \cdots & 
\frac{Q_m\hat{X}_{t-1}}{{\bf 1}_n^\top Q_m\hat{X}_{t-1}}
\end{array}
\right]Y_t \label{eq:genfilter}
\end{equation}
with the initial condition determined by the initial law of $Z$. The 
prediction $\hat{X}_{t+1|t}$ is equal to $A\hat{X}_t$ with $A=\sum_{i=1}^m 
Q_i$ and $X_{0|-1}=EX_0=p_0$. For the prediction $\hat{Y}_{t+1|t}$ 
we have $\hat{Y}_{t+1|t}=C\hat{X}_t$ with $C=(I_m\otimes {\bf 1}_n^\top )\bar{Q}$. 
\label{thm:filterp}
\end{thm}
{\bf Proof.}
We use equation  (\ref{eq:bayes})   with ${\cal F}^0 = {\cal 
F}^Y_t$, ${\cal  
H}=\sigma(Y_{t+1})$, which is generated by the sets $H_i=\{Y_{t+1}=f_i\}$ and
$U=X_{t+1}$. Thus we obtain
\[
E[X_{t+1}|{\cal F}^Y_{t+1}] = \sum_{i=1}^m E_i[X_{t+1}|{\cal F}^Y_t]1_{H_i} 
= \sum_{i=1}^m \frac{E[X_{t+1}1_{H_i}|{\cal F}^Y_t]}{P(H_i|{\cal 
F}^Y_t)}1_{H_i}.
\] Then we use the Markov property of $Z$ to write 
\[
E[X_{t+1}1_{H_i}|{\cal F}^Y_t] = E[E[X_{t+1}1_{H_i}|{\cal F}_t]|{\cal F}^Y_t] =
E[Q_iX_t|{\cal F}^Y_t] = Q_i\hat{X}_t. \]
Since $P(H_i|{\cal F}^Y_t) = E[1_{H_i}|{\cal F}_t]={\bf 
1}_n^\top E[X_{t+1}1_{H_i}|{\cal F}^Y_t]$ we get equation (\ref{eq:genfilter}).
\\
Define now $A=\sum_{i=1}^m Q_i=({\bf 1}_m^\top \otimes I_n)\bar{Q}$ and 
$C=(I_m\otimes {\bf 
1}_n^\top )\bar{Q}$. Then we have $E[X_{t+1}  |{\cal F}_t]=AX_t$ and 
$E[Y_{t+1}|{\cal F}_t]= CX_t$. As a consequence we get by reconditioning that
$\hat{X}_{t+1|t}=A\hat{X}_t$ and that $\hat{Y}_{t+1|t}=C\hat{X}_t$. \hfill 
$\square$ 
\bigskip \\
We see that the filter $\hat{X}_t$ satisfies a completely recursive system, 
that is, $\hat{X}_t$ is completely determined by $\hat{X}_{t-1}$ and $Y_t$.
In absence of further conditions on the matrix $Q$ (in particular the 
factorization property) there seems to be no complete recursion that is 
satisfied by $X_{t|t-1}$. The reason for this is that we don't have the Markov 
property of $W$ with respect to $\mathbb{G}$, unless the factorization property 
holds,  
in which case the formulas above take a particular nice form. See  
subsection \ref{subsection:filterhmc}. 
\begin{remark}\label{remark:altfilter}
\textup{
It follows from equation~(\ref{eq:intQ}) that the filter (\ref{eq:genfilter}) can 
alternatively be expressed as
\[
\hat{X}_t=
\left[
Q_1\hat{X}_{t-1},\ldots, 
Q_m\hat{X}_{t-1}
\right]
\diag(C\hat{X}_{t-1})^{-1}Y_t.
\]
Indeed, from equation~(\ref{eq:intQ}) we obtain ${\bf 1}^\top Q_i X_t= 
P(Y_{t+1}=f_i|\mathcal{F}_t)$, hence $E[Y_{t+1}|\mathcal{F}_t]$ is the vector with 
elements ${\bf 1}^\top Q_i X_t$. Conditioning of this vector on $\mathcal{F}^Y_t$ 
gives that $\hat{Y}_{t+1|t}$ is the vector with elements  ${\bf 1}^\top Q_i 
\hat{X}_t$. So we can rewrite~(\ref{eq:genfilter}) as  $\hat{X}_t=
\left[
Q_1\hat{X}_{t-1},\ldots,
Q_m\hat{X}_{t-1}\right]
\diag(Y_{t|t-1})^{-1}Y_t$
and the result follows.
}
\end{remark}

\subsection{Filter for $\Sigma_S$} \label{subsection:filters}
In this section we obtain the filter for a system in $\Sigma_S$, so we work 
with a Markov chain $W_t=X_t\otimes Y_{t-1}$ with transition matrix 
$R=\bar{R}({\bf 1}_m^\top \otimes I_n)$, where 
the matrix $\bar{R}$ can be written as
\begin{equation}
\bar{R}=
\left[
\begin{array}{c}
R_1 \\
\vdots \\
R_m
\end{array}
\right].
\end{equation}
for certain matrices  $R_i$  in $\mathbb{R}^{n\times n}$. No further assumptions 
on the $R_i$ are made. Observe that the $R_i$ have the interpretation that 
\begin{equation}\label{eq:intR}
R_iX_t =  E[X_{t+1}1_{\{Y_{t}=f_i\}}|{\cal G}_t].
\end{equation}
Then we have
\begin{thm}
The predictor $\hat{X}_{t|t-1}$ is given by the recursion
\begin{equation}
\hat{X}_{t+1|t}=
\left[
\begin{array}{ccc}
\frac{R_1\hat{X}_{t|t-1}}{{\bf 1}_n^\top R_1\hat{X}_{t|t-1}} & \cdots & 
\frac{R_m\hat{X}_{t|t-1}}{{\bf 1}_n^\top R_m\hat{X}_{t|t-1}}
\end{array}
\right]Y_t \label{eq:genfilters}
\end{equation}
with the initial condition $X_{0|-1}=EX_0$. 
For the filter $\hat{X}_t$ and for $\hat{Y}_{t+1|t}$ we have the following relations. 
\begin{equation}
\hat{X}_{t+1}=\diag(\hat{X}_{t+1|t})G^\top \diag(\hat{Y}_{t+1|t})^{-1}Y_{t+1}, 
\label{eq:filters}
\end{equation}
where $G=(I_m\otimes {\bf 1}_n^\top )\bar{R}$ and
\begin{equation}
\hat{Y}_{t+1|t}=GX_{t+1|t}. \label{eq:predy}
\end{equation}
\label{thm:filters}
\end{thm}
{\bf Proof.} 
We use equation (\ref{eq:bayes}) with ${\cal F}^0 = {\cal 
F}^Y_{t-1}$, ${\cal  
H}=\sigma(Y_{t})$, which is generated by the sets $H_i=\{Y_{t}=f_i\}$ and
$U=X_{t+1}$. Then we obtain
$E[X_{t+1}|{\cal F}^Y_{t}] = \sum_{i=1}^m E_i[X_{t+1}|{\cal F}^Y_{t-1}]1_{H_i} 
= \sum_{i=1}^m \frac{E[X_{t+1}1_{H_i}|{\cal F}^Y_{t-1}]}{P(H_i|{\cal 
F}^Y_t)}1_{H_i}$. Then we use the Markov property of $W$ to write 
\begin{eqnarray*}
E[X_{t+1}1_{H_i}|{\cal F}^Y_{t-1}] & = & E[E[X_{t+1}1_{H_i}|{\cal G}_t]|{\cal 
F}^Y_{t-1}] \\
& =& 
E[R_iX_t|{\cal F}^Y_{t-1}] \\
& = & R_i\hat{X}_{t|t-1}.
\end{eqnarray*}
Since $P(H_i|{\cal F}^Y_{t-1}) = E[1_{H_i}|{\cal F}^Y_{t-1}]={\bf 
1}_n^\top E[X_{t+1}1_{H_i}|{\cal F}^Y_{t-1}]$ we get equation 
(\ref{eq:genfilters}). \\
To derive the formula (\ref{eq:filters}) for the filter we proceed similarly, 
using lemma  
\ref{lemma:bayes} again with $U=X_{t+1}$, ${\cal F}^0={\cal F}^Y_t$ and 
${\cal H}=\sigma(Y_{t+1})$ generated by the sets $H_i=\{Y_{t+1}=f_i\}$. Then 
we can write equation (\ref{eq:bayes}) as $E[X_{t+1}|{\cal F}^Y_{t+1}] = 
E[X_{t+1}Y^\top _{t+1}|{\cal F}^Y_t]\diag(\hat{Y}_{t+1|t})^{-1}Y_{t+1}$. 
\begin{eqnarray*}
E[X_{t+1}Y^\top _{t+1}|{\cal F}^Y_t] & = & E[E[X_{t+1}Y^\top _{t+1}|{\cal 
G}_{t+1}]|{\cal F}^Y_t] \\
& = & E[X_{t+1}E[Y^\top _{t+1}|{\cal 
G}_{t+1}]|{\cal F}^Y_t] \\
& = & E[X_{t+1}(GX_{t+1})^\top |{\cal F}^Y_t]  \\
& = & E[\diag(X_{t+1})|{\cal F}^Y_t]G^\top .
\end{eqnarray*}
Then equation (\ref{eq:filters}) follows, 
as well as equation (\ref{eq:predy}), since we have $E[Y_{t+1}|{\cal 
F}^Y_t] =
E[Y_{t+1}X_{t+1}^\top |{\cal F}^Y_t]{\bf 1}_n = G\diag(\hat{X}_{t+1|t}){\bf 
1}_n  = G\hat{X}_{t+1|t}$. 
\hfill $\square$
\begin{remark}
\textup{
By a similar argument as in remark~\ref{remark:altfilter} we can rewrite the 
recursion~(\ref{eq:genfilters}) for the predictor as 
\[
\hat{X}_{t+1|t}=
\left[
R_1\hat{X}_{t|t-1},\dots,R_m\hat{X}_{t|t-1}\right]\diag(G\hat{X}_{t|t-1})^{-1}Y_t.
\]
}
\end{remark}
\begin{remark}
\textup{
Notice that in contrast with what we got in subsection 
\ref{subsection:filterp} for 
$\Sigma_P$ here the {\em predictor} satisfies a completely recursive system, 
whereas we obtain the filter in terms of the predictor.
}
\end{remark}
The formulas above take a particular nice form if the system satisfies the 
splitting property. See subsection \ref{subsection:filterhmc}.

\subsection{Filter for a Hidden Markov Chain} \label{subsection:filterhmc}
In this section we return to the setting of sections 
\ref{section:preliminaries}  
and \ref{section:alternative} and we give the 
recursive filtering formula for the stochastic  
system with the HMC $Y$ as its output. Therefore, we can apply the results of 
subsection \ref{subsection:filterp} with the specification that 
$\bar{Q}=\Delta(G)A$, so we have $Q_i = \diag(G_{i.})A$ and ${\bf 1}_n^\top Q_i = 
G_{i.}A$. 
The following holds.
\begin{thm}
(i) The conditional distribution of the $X_t$ given $Y_0,\ldots,Y_t$ is 
recursively determined by                                                
\begin{equation}
\hat{X}_{t}=\mbox{{\em diag}}(A\hat{X}_{t-1})G^\top \mbox{{\em diag}}(GA\hat{X}_{t-1})^{-
1}Y_t, \label{eq:filter}
\end{equation}
with initial condition $\hat{X}_0=\mbox{{\em diag}}(p_0) G^\top \mbox{{\em diag}}(Gp_0)^{-1}Y_0$, with 
$p_0=EX_0$. \\
(ii) The conditional distribution of the $X_t$ given $Y_0,\ldots,Y_{t-1}$ is 
recursively determined by                                                
\begin{equation}
\hat{X}_{t+1|t}=A\mbox{{\em diag}}(\hat{X}_{t|t-1})G^\top \mbox{{\em 
diag}}(G\hat{X}_{t|t-1})^{-
1}Y_t,    \label{eq:predhmc}
\end{equation}
with initial condition $X_{0|-1}=EX_0=p_0$. \\
(iii) The conditional expectation $\hat{Y}_{t+1|t} = E[Y_{t+1}|{\cal F}^Y_t]$ 
is given by 
\begin{equation}
\hat{Y}_{t+1|t} = GA\mbox{{\em diag}}(\hat{X}_{t|t-1})G^\top \mbox{{\em 
diag}}(G\hat{X}_{t|t-
1})^{-1}Y_t             \label{eq:ycondexp}
\end{equation}
\label{thm:filterhmc}
\end{thm}
{\bf Proof}. (i) Just use equation (\ref{eq:genfilter}) and notice that 
\[ \frac{Q_i\hat{X}_{t-1}}{{\bf 1}_n^\top Q_i\hat{X}_{t-1}} =\diag(A\hat{X}_{t-
1})\frac{G_{i.}^\top }{(GA\hat{X}_{t-1})_i}. \]
(ii) follows from (i), since we know from theorem \ref{thm:filterp} that 
$\hat{X}_{t+1|t}=A\hat{X}_t$. \\ 
(iii) also follows from theorem \ref{thm:filterp}, upon noticing that $C$ 
now becomes $GA$ in view of (\ref{eq:deltarule}). \hfill $\square$ 
\begin{remark}
\textup{
Here both the filter and the predictor satisfy a 
complete  recursive system. This is not surprising, because a HMC is a stochastic 
system belonging to both $\Sigma_P$ and $\Sigma_S$. Notice that
theorem 
\ref{thm:filterhmc} can alternatively be derived from theorem \ref{thm:filters}, since under 
the assumptions of the present subsection we have that $R_i=A\diag(G_{i.})$. 
}
\end{remark}
\begin{remark}
\textup{
If we define for $x \in \mathbb{R}^n_+$ the matrix 
\[
G_x:= 
\diag(x)G^\top \diag(Gx)^{-1},
\] 
then equations (\ref{eq:filter}),  
(\ref{eq:predhmc}) and  (\ref{eq:ycondexp})  take 
the form $\hat{X}_t=G_{A\hat{X}_{t-1}}Y_t$, $\hat{X}_{t+1|t}= AG_{\hat{X}_{t|t-1}}Y_t$ 
and $\hat{Y}_{t+1|t}= GAG_{\hat{X}_{t|t-1}}Y_t$. \\
One may check that under the condition that $Y$ is a {\em 
deterministic} function of $X$ (in which case the columns of $G$ are 
basis vectors of $\mathbb{R}^m$) the matrices $G_x$ 
are 
right pseudo-inverses of $G$.
}
\end{remark}

\appendix

\section{A lemma on conditional expectations}

Consider some probability space $(\Omega,\mathcal{F},P)$ and
let ${\cal H}$ be a sub-$\sigma$-algebra of ${\cal F}$ that is generated by 
a finite partition $\{H_1,\ldots,H_k\}$ of $\Omega$, satisfying $P(H_i)>0$ for 
all $i$. We introduce the (conditional)
probability measures $P_i$ on $(\Omega,{\cal F})$ defined by $P_i(F)= E[1_F 
\frac{1_{H_i}}{P(H_i)}] = P(F|H_i)$. Expectation with respect to $P_i$ is 
denoted by $E_i$. Notice that for a random variable $X$ with finite expectation we have
\begin{equation}\label{eq:expi}
E1_{H_i}X = P(H_i)E_iX.
\end{equation}
We also have that for any sub-$\sigma$-algebra $\mathcal{F}_0$ and an integrable 
random variable $X$ the equality
\begin{equation}\label{eq:condi}
E[X1_{H_i}|\mathcal{F}_0]1_{H_i} = P(H_i|\mathcal{F}_0)E_i[X|\mathcal{F}_0].
\end{equation}
Recall that for any integrable random 
variable $U$  it holds that
\begin{equation}
E[U|{\cal H}]= \sum_{i=1}^kE_i[U]1_{H_i}.       \label{eq:conex}
\end{equation}
We extend this result in the following easy to prove lemma. It is used frequently in 
sections~\ref{section:alternative}, \ref{section:stochsyst} and~\ref{section:filter}.
\begin{lemma}\label{lemma:bayes} 
Let $\mathcal{F}^0$ be some sub-$\sigma$-algebra of $\mathcal{F}$. Then the following equalities hold true.
\begin{eqnarray}
E[U|{\cal F}^0\vee {\cal H}] & = & \sum_{i=1}^k E_i[U|{\cal 
F}^0]1_{H_i}.  \label{eq:bayes}  \\
E[U|{\cal F}^0] & = & \sum_{i=1}^k E_i[U|{\cal 
F}^0]E[1_{H_i}|{\cal F}^0]. \label{eq:bayes3} \\                      
E[1_{H_j}U|{\cal F}^0] & = & E[1_{H_j}|{\cal F}^0] E_j[U|{\cal F}^0]. 
\label{eq:bayes2} 
\end{eqnarray}
\end{lemma}
{\bf Proof.}  Concerning the first equality we have to show that 
\[
 E\{1_{F\cap H_j} U\} = E\{1_{F\cap 
H_j}\sum_{i=1}^k E_i[U|{\cal F}^0]1_{H_i}\} 
\] for all 
$F\in {\cal F}^0$, because every set in $\mathcal{F}^0\vee\mathcal{H}$ can be written as 
a finite union of sets $F\cap H_j$ with some $F\in\mathcal{F}^0$ and because the RHS of~(\ref{eq:bayes}) is
clearly ${\cal F}^0\vee {\cal H}$-measurable.
We develop
\begin{eqnarray*}
\mbox{$E\{1_{F\cap H_j}\sum_{i=1}^k E_i[U|{\cal F}^0]1_{H_i}\}$} 
& = & E\{1_{F\cap H_j}E_j[U|{\cal F}^0]\} \\
& = & E_j\{1_FE_j[U|\mathcal{F}^0]\}P(H_j) \\
& = & E_j\{1_F U\}P(H_j) \\
& = & E\{1_{F\cap H_j}U\}.
\end{eqnarray*}
In these computations we used~(\ref{eq:expi}) in the second and fourth equality and the 
defining property of conditional expectation in  the third.
This proves (\ref{eq:bayes}). \\
The second equality is a direct consequence of the first by conditioning on $\mathcal{F}^0$.
The third equality follows from the second one
by taking $1_{H_j}U$ instead 
of $U$.
\hfill $\square$
\begin{remark}\label{remark1}
\textup{ 
If we take in lemma \ref{lemma:bayes} ${\cal F}^0$ the trivial 
$\sigma$-algebra, then (\ref{eq:bayes}) reduces to (\ref{eq:conex}). If 
$P(H_i)=0$ for some $i$, then $P_i$ is not well defined but (\ref{eq:bayes}) 
is still valid provided we  
define $E_i[U|{\cal F}^0]$ to be zero for such an $i$. 
}
\end{remark}
\begin{remark}
\textup{
Equation (\ref{eq:bayes2}) is also known as the conditional Bayes theorem, 
cf. \cite[page 23]{elliottetal}.
}
\end{remark}


\begin{thebibliography}{99}
\bibitem{baumpetrie}
L.E.\ Baum and T.\ Petrie (1966), Statistical inference for probabilistic 
functions of finite state Markov chains, {\em Ann.\ Math.\ Stat.} 37, pp. 1554-
1563.
\bibitem{elliottetal}
R.\ J.\ Elliott, L.\ Aggoun and J.\ B.\ Moore (1995), {\em Hidden Markov 
Models. Estimation and Control}, Springer.
\bibitem{finesso}
L.\ Finesso (1990), {\em Consistent Estimation of the Order for Markov and 
Hidden Markov Chains}, Dissertation University of Maryland.
\bibitem{magnusneudecker}
J.R.\ Magnus \& H.\ Neudecker (1988), {\em Matrix Differential Calculus with 
Applications in Statistics and Econometrics}, Wiley.
\bibitem{picci}
G.\ Picci (1978),  On the internal structure of finite state stochastic 
processes, in {\em Recent Developments in Variable Structure Systems}, 
Springer Lecture Notes in Economics and Math. systems, Vol. 162.
\bibitem{schuppen}
J.H.\ van Schuppen (1989), Stochastic realization problems, In {\em 
Three decades of Mathematical System Theory}, H. Nijmeijer, J.M. Schumacher 
(Eds.), Springer Lecture Notes in Control and Information Sciences 135, pp. 
480-523. 
\bibitem{spreij}
P.J.C.\ Spreij (2001), On the Markov property of a hidden Markov chain, {\em 
Statistics and Probability Letters}, Vol 52/3,   pp 279-288.
\end{thebibliography}
\end{document}